\documentclass[12pt]{article}
\usepackage{amsmath}
\usepackage{prodint}
\usepackage{graphicx}
\usepackage{enumerate}
\usepackage{natbib}
\usepackage{url} 
\usepackage{binomexp}
\usepackage{mathtools}
\usepackage{amsfonts}
\usepackage{mathrsfs}
\RequirePackage{subcaption}
\usepackage{booktabs}
\usepackage{bbm}
\usepackage{dsfont}
\usepackage{longtable}
\usepackage{tabulary}
\usepackage{yhmath}
\usepackage{bm}
\usepackage{amssymb}
\usepackage[nodisplayskipstretch]{setspace}
\usepackage{tikz}
\usetikzlibrary{arrows.meta, decorations.pathmorphing, decorations.pathreplacing, positioning, calc, shadows.blur}
\usepackage{float}
\usepackage{color}
\usepackage{stmaryrd}
\usepackage{paralist}
\usepackage{multirow}
\usepackage{threeparttable}
\usepackage{soul}
\usepackage{algorithm,algpseudocode}
\usepackage{amsthm}
\usepackage{array}
\DeclareFontFamily{U}{mathx}{\hyphenchar\font45}
\DeclareFontShape{U}{mathx}{m}{n}{<5> <6> <7> <8> <9> <10> <10.95> <12> <14.4> <17.28> <20.74> <24.88> mathx10}{}
\DeclareSymbolFont{mathx}{U}{mathx}{m}{n}
\DeclareMathAccent{\widecheck}{0}{mathx}{"71}

\numberwithin{equation}{section}
\theoremstyle{plain}
\newtheorem{theorem}{Theorem}[section]
\newtheorem{lemma}[theorem]{Lemma}
\newtheorem{proposition}[theorem]{Proposition}

\theoremstyle{definition}
\newtheorem{assumption}{Assumption}
\newtheorem{condition}{Condition}
\newtheorem{definition}[theorem]{Definition}

\newcommand{\R}{\mathbb{R}}
\newcommand{\E}{\mathbb{E}}
\newcommand{\pr}{\mathbb{P}}
\newcommand{\1}{\mathds{1}}
\newcommand{\wto}{\rightsquigarrow}
\newcommand{\bX}{\bm{X}}
\newcommand{\be}{\bm{e}}
\newcommand{\bw}{\bm{w}}
\newcommand{\bg}{\bm{g}}
\newcommand{\bv}{\bm{v}}
\newcommand{\bbeta}{\bm{\beta}}
\newcommand{\bpsi}{\bm{\psi}}
\newcommand{\bvarepsilon}{\bm{\varepsilon}}
\newcommand{\bsk}{\bm{\mathsf{k}}}
\newcommand{\bso}{\bm{\mathsf{o}}}
\newcommand{\bzero}{\bm{0}}
\newcommand{\mG}{\mathbf{G}}
\newcommand{\mSig}{\bm{\Sigma}}
\newcommand{\bbh}{\widehat{\bbeta}}
\newcommand{\that}{\widehat{\tau}}
\newcommand{\Vh}{\widehat{V}}
\newcommand{\cstar}{c_{\ast}}
\newcommand{\Phibar}{\overline{\Phi}}
\newcommand{\var}{\operatorname{var}}

\usepackage[colorlinks=true,linkcolor=blue,citecolor=black]{hyperref}

\newcommand{\blind}{0}

\allowdisplaybreaks

\begin{document}

\def\spacingset#1{\renewcommand{\baselinestretch}%
{#1}\small\normalsize} \spacingset{1}


\if0\blind
{
  \begin{center}
	{\Large \bf Studentized inference in randomized experiments when a few units dominate the outcome variation}\\[0.5cm]
	{\large Xinyuan Chen$^{\ast}$ ~ Jonathan W. Kiesel ~ Jonathan R. Woody}\\[0.2cm]
	Department of Mathematics and Statistics, Mississippi State University, \\Mississippi State, MS, USA\\
	$^\ast$xchen@math.msstate.edu
  \end{center}
} \fi

\if1\blind
{
  \bigskip
  \bigskip
  \bigskip
  \begin{center}
    {\Large \bf Studentized inference in randomized experiments when a few units dominate the outcome variation}
\end{center}
  \medskip
} \fi

\begin{abstract}
    In completely randomized experiments with heavy-tailed outcomes, a few units can carry a share of the residual variation that does not vanish as the number of units grows. In this setting, we study the studentized difference in means and Lin's regression-adjusted estimator under a design-based framework. The studentized statistic has a nonnormal limit in general.  Populations with the same marginal distributions of potential outcomes can have different limits, as the limit depends on how the two potential outcomes of each dominant unit are paired. When the outcome errors are drawn from a heavy-tailed law, the limit is random and, outside a degenerate case, varies with the population. Under balanced allocation, the conventional Wald interval has an explicit asymptotic coverage floor. At the nominal level of $0.95$, the floor is about $0.42$ for every population and about $0.84$ when the two potential residuals of each unit are equal. Both floors are close to the best possible, and without balanced allocation, the coverage can fall to the smaller of the proportions of treated and control units.
\end{abstract}

\noindent%

\spacingset{1.5}

\section{Introduction}
\label{sec:intro}

Randomized experiments are commonly analyzed with the difference in means or with a regression-adjusted estimator, together with a Wald interval based on a Neyman-type variance estimator. In the design-based framework, the covariates and potential outcomes are fixed, the treatment assignment is the only source of randomness, and the interval is justified by a finite-population central limit theorem \citep{Lin2013, LiDing2017}. Such theorems require that no unit carry a share of the outcome variation that stays bounded away from zero as the number of units grows, a condition known as negligibility \citep{Noether1949, ErdosRenyi1959, Hajek1960}. In practice, however, outcomes such as revenue, spending, and charges are heavy-tailed, and a few units can carry much of the variation. As \citet{Young2019} noted, in published experiments, deleting a single observation can make about a third of the results that are significant at the level $0.01$ insignificant. For example, in the RAND Health Insurance Experiment \citep{AronDineEinavFinkelstein2013} studied in Section \ref{sec:data}, one person carries about a third of the residual variation of second-year medical spending.

Without negligibility, sums over a randomly selected subset of a finite population can have nonnormal limits with discrete components \citep{Rosen1965, Kallenberg1973, Jirina1982, BerkesHorvathSchauer2010}, and permuted sums of heavy-tailed samples have limits that depend on the sample \citep{AueBerkesHorvath2008}. For independent observations without finite variance, the limits of self-normalized sums are known \citep{Logan1973}. Student's $t$ statistic is asymptotically normal if and only if the law of the observations is in the domain of attraction of the normal law \citep{GineGotzeMason1997}, and under symmetry, it has explicit tail bounds \citep{Efron1969, Edelman1990, Pinelis1994}. The Berry--Esseen bound of \citet{Bloznelis1999} for the finite-population Student statistic implies asymptotic normality under the condition of \citet{ErdosRenyi1959}, and the theory of studentized permutation and randomization tests \citep{Janssen1997, ChungRomano2013, WuDing2021} establishes validity when the studentized statistic is asymptotically normal. Methods designed for heavy-tailed outcomes change either the estimator or the estimand \citep{Athey2023, SuLi2024, Ghosh2026}. Concentration inequalities give intervals with finite-sample coverage \citep{Freidling2026}, but they require known bounds on the potential outcomes. None of the above results gives the randomization law of the conventional studentized statistic when a few units dominate the residual variation. The limits for random subsets concern sums of one value per unit, whereas the studentized statistic is a ratio whose numerator and denominator depend on both potential outcomes of the dominant units.

This article studies the studentized difference in means and the interacted estimator of \citet{Lin2013}, with covariates whose leverage vanishes, when the residual variation is not negligible. The studentized statistic has a limit law determined by the limits of the normalized residual pairs of the dominant units and of the second-moment matrix of the remaining units, provided the limit of the normalized variance estimator is positive with probability one (Theorem \ref{thm:studentized}). The law is nonnormal in general. It depends on how the two potential residuals of each dominant unit are paired, not only on their marginal distributions, so populations with the same marginal distributions of potential outcomes can have different limiting coverage (Proposition \ref{prop:pairing}). When the outcome errors are drawn from a law with tail index between one and two, the randomization law of the statistic for the population at hand converges to a random law determined by the largest residual pairs of the realized population, and this law varies with the population unless the spectral measure of the errors is concentrated on a line determined by the allocation (Proposition \ref{prop:sampling}). Under balanced allocation, the Wald interval has an explicit asymptotic coverage floor for every population, about $0.42$ at the nominal level of $0.95$ (Proposition \ref{prop:balancedfloor}), and a higher floor, about $0.84$, when the two potential residuals of each unit are equal (Theorem \ref{thm:floor}). Both floors are close to the best possible, and neither holds at the nominal level of $0.95$ when $40\%$ of the units are treated. Without balanced allocation, the coverage can fall to the smaller of the proportions of treated and control units (Proposition \ref{prop:counter}). The floors are randomization analogs of the tail bounds for Student's $t$ statistic under symmetry, with the balance of the design in place of symmetric errors.

Simulations with heavy-tailed errors (Section \ref{sec:sim}) and a data application (Section \ref{sec:data}) show how closely the limits describe populations of moderate size and how the coverage of the Wald interval behaves when a few units dominate.

\section{Setup}
\label{sec:setup}

\subsection{Finite populations and the design}
\label{sec:finitepop}

For each $N\ge2$, let $\Pi_N=\{(\bX_i,Y_i(1),Y_i(0)):i\le N\}$ be a finite population of covariates $\bX_i\in\R^d$ and potential outcomes $Y_i(1),Y_i(0)\in\R$. These quantities are fixed but may change with $N$. We suppress their dependence on $N$ in the notation. The dimension $d\ge0$ is fixed, and $d=0$ means that there are no covariates. The estimand is the sample average treatment effect,
\begin{equation*}
    \tau_N=\frac1N\sum_{i=1}^N\{Y_i(1)-Y_i(0)\}.
\end{equation*}

\begin{assumption}[Complete randomization]
    \label{ass:design}
    For each $N$, $1\le N_1\le N-1$, and the vector of treatment indicators $(Z_1,\dots,Z_N)$ is uniformly distributed over all binary vectors of length $N$ with $N_1$ entries equal to one. Also, $N_0=N-N_1$, and $\pi_N=N_1/N\to\pi\in(0,1)$.
\end{assumption}

\subsection{Projections and leverage}
\label{sec:projections}

Let $\overline{Y}(z)=N^{-1}\sum_iY_i(z)$ for $z\in\{0,1\}$. If $d\ge1$, let $\overline{\bX}=N^{-1}\sum_i\bX_i$ and $\mG=\sum_i(\bX_i-\overline{\bX})(\bX_i-\overline{\bX})^\top$, which we assume to be nonsingular for every $N$, and define the finite-population projection coefficients and residuals by
\begin{equation}
    \label{eq:proj}
    \bbeta_z=\mG^{-1}\sum_{i=1}^N(\bX_i-\overline{\bX})Y_i(z),
    \qquad
    e_i(z)=Y_i(z)-\overline{Y}(z)-(\bX_i-\overline{\bX})^\top\bbeta_z .
\end{equation}
If $d=0$, let $e_i(z)=Y_i(z)-\overline{Y}(z)$. No outcome model is assumed. By construction,
\begin{equation}
    \label{eq:orth}
    \sum_{i=1}^Ne_i(z)=0,
    \qquad
    \sum_{i=1}^N(\bX_i-\overline{\bX})e_i(z)=\bzero,
\end{equation}
where the second identity applies if $d\ge1$. The leverage of unit $i$ is $h_i=(\bX_i-\overline{\bX})^\top\mG^{-1}(\bX_i-\overline{\bX})$ if $d\ge1$, and $h_i=0$ if $d=0$. Then, $0\le h_i\le1$ and $\sum_ih_i=d$. Write $\be_i=(e_i(1),e_i(0))^\top$ for the residual pair of unit $i$, and let
\begin{equation}
    \label{eq:scales}
    r_N^2=\sum_{i=1}^N\|\be_i\|^2,
    \qquad
    c_i=\frac{e_i(1)}{\pi_N}+\frac{e_i(0)}{1-\pi_N}.
\end{equation}
The quantity $r_N^2$ is the total residual variation, and $c_i$ combines the two potential residuals of unit $i$.

\subsection{The estimator and the Wald interval}
\label{sec:estimators}

For $z\in\{0,1\}$, let $\overline{\bX}_z$ and $\overline{Y}_z$ denote the averages of the covariates and of the observed outcomes $Y_i=Z_iY_i(1)+(1-Z_i)Y_i(0)$ over $\{i:Z_i=z\}$. If $d\ge1$, define the within-arm Gram matrix and slope by
\begin{equation}
    \label{eq:gram}
    \mG_z=\sum_{i:Z_i=z}(\bX_i-\overline{\bX}_z)(\bX_i-\overline{\bX}_z)^\top,
    \qquad
    \bbh_z=\mG_z^{-1}\sum_{i:Z_i=z}(\bX_i-\overline{\bX}_z)Y_i,
\end{equation}
with $\bbh_z=\bzero$ if $\mG_z$ is singular. The estimator is
\begin{equation}
    \label{eq:tau}
    \that=\overline{Y}_1-\overline{Y}_0-(\overline{\bX}_1-\overline{\bX})^\top\bbh_1+(\overline{\bX}_0-\overline{\bX})^\top\bbh_0
\end{equation}
if $d\ge1$, and the difference in means $\that=\overline{Y}_1-\overline{Y}_0$ if $d=0$. If $\mG_1$ and $\mG_0$ are nonsingular, $\that$ is the coefficient of $Z_i$ in the least squares regression of $Y_i$ on $(1,Z_i,\bX_i-\overline{\bX},Z_i(\bX_i-\overline{\bX}))$, the interacted estimator of \citet{Lin2013}. Suppose $N_1,N_0\ge2$. For $Z_i=z$, let $\widehat e_i=Y_i-\overline{Y}_z-(\bX_i-\overline{\bX}_z)^\top\bbh_z$ if $d\ge1$ and $\widehat e_i=Y_i-\overline{Y}_z$ if $d=0$. Define
\begin{equation}
    \label{eq:Vhat}
    \Vh=\frac{s_1^2}{N_1}+\frac{s_0^2}{N_0},
    \qquad
    s_z^2=\frac{1}{N_z-1}\sum_{i:Z_i=z}\widehat e_i^{\,2},
    \qquad
    T_N=\frac{\that-\tau_N}{\Vh^{1/2}}.
\end{equation}
If $d=0$, $\Vh$ is the Neyman variance estimator. If $\Vh=0$, we set $T_N=0$ when $\that=\tau_N$, and $T_N=+\infty$ or $T_N=-\infty$ according to the sign of $\that-\tau_N$ otherwise, so that $\{|T_N|\le t\}$ is the event that $\that\pm t\Vh^{1/2}$ covers $\tau_N$, for $t>0$. The Wald interval of nominal level $1-\varsigma$ is $\that\pm\Phi^{-1}(1-\varsigma/2)\Vh^{1/2}$, where $\Phi$ is the standard normal distribution function.

\subsection{Exact decompositions}
\label{sec:decomp}

Define the residual contrast and, if $d\ge1$, the remainder from estimating the slopes by
\begin{equation}
    \label{eq:RN}
    L_N=\sum_{i=1}^N(Z_i-\pi_N)\left\{\frac{e_i(1)}{N_1}+\frac{e_i(0)}{N_0}\right\},
    \qquad
    R_N=-(\overline{\bX}_1-\overline{\bX})^\top(\bbh_1-\bbeta_1)+(\overline{\bX}_0-\overline{\bX})^\top(\bbh_0-\bbeta_0),
\end{equation}
and let $R_N=0$ if $d=0$. Let $\overline{e}_z(z)=N_z^{-1}\sum_{i:Z_i=z}e_i(z)$. If $d=0$, let $\mathsf{C}_{N,z}=N_z\overline{e}_z(z)^2$. If $d\ge1$, let $\bsk_z=\sum_{i:Z_i=z}(\bX_i-\overline{\bX}_z)e_i(z)$ and
\begin{equation}
    \label{eq:Cz}
    \mathsf{C}_{N,z}=N_z\overline{e}_z(z)^2+\bsk_z^\top\mG_z^{-1}\bsk_z
\end{equation}
on the event that $\mG_z$ is nonsingular, and $\mathsf{C}_{N,z}=0$ otherwise.

\begin{proposition}
    \label{prop:decomp}
    Suppose Assumption \ref{ass:design} holds. For every assignment,
    \begin{equation}
        \label{eq:decint}
        \that-\tau_N=L_N+R_N,
        \qquad
        NL_N=\sum_{i=1}^N(Z_i-\pi_N)c_i,
    \end{equation}
    and $\E(NL_N)=0$ and $\var(NL_N)=\{N_1N_0/(N(N-1))\}\sum_ic_i^2$. If $N_1,N_0\ge2$ and either $d=0$ or $\mG_1$ and $\mG_0$ are nonsingular, then $\mathsf{C}_{N,1},\mathsf{C}_{N,0}\ge0$ and
    \begin{equation}
        \label{eq:NVhat}
        N^2\Vh=\sum_{z\in\{0,1\}}\frac{N^2}{N_z(N_z-1)}\left\{\sum_{i:Z_i=z}e_i(z)^2-\mathsf{C}_{N,z}\right\}.
    \end{equation}
\end{proposition}

Section S2 proves Proposition \ref{prop:decomp}. The estimation error splits into the residual contrast $L_N$, which is the difference in means computed from the population residuals, and the remainder $R_N$ from estimating the slopes. In each arm, the variance estimator sums the squared population residuals of the observed units, minus $\mathsf{C}_{N,z}$, the part of that sum explained by regressing these residuals on an intercept and the covariates within arm $z$. Both $L_N$ and the leading term of \eqref{eq:NVhat} depend on the population only through the residual pairs $\be_i$.

\begin{condition}[Vanishing leverage]
    \label{cond:lev}
    The leverages satisfy $\max_ih_i\to0$.
\end{condition}

Condition \ref{cond:lev} holds if $d=0$. If $d\ge1$, it holds if $N^{-1}\max_i\|\bX_i-\overline{\bX}\|^2\to0$ and $\liminf_N\lambda_{\min}(N^{-1}\mG)>0$, because $h_i\le\|\bX_i-\overline{\bX}\|^2/\lambda_{\min}(\mG)$. It restricts the covariates only.

\begin{lemma}
    \label{lem:remainder}
    Suppose Assumption \ref{ass:design} and Condition \ref{cond:lev} hold. Then, if $d\ge1$, the probability that $\mG_1$ or $\mG_0$ is singular tends to zero. For every $d\ge0$, if $r_N>0$ for every $N$, then $NR_N/r_N\to0$ and $\mathsf{C}_{N,z}/r_N^2\to0$ in probability, for $z\in\{0,1\}$.
\end{lemma}

Section S2.3 proves Lemma \ref{lem:remainder} with a bound for each $N$. Under Condition \ref{cond:lev}, estimating the slopes changes the numerator and the denominator of $T_N$ only at a smaller order than the residual scale, regardless of the outcomes. The limits below are therefore determined by the residual pairs $\be_i$.

\section{The limit of the studentized statistic}
\label{sec:studentized}

\subsection{Residual configurations}
\label{sec:config}

When a few units carry a share of $r_N^2$ that does not vanish, the numerator $N(\that-\tau_N)$ and the denominator $N^2\Vh$ of $T_N$ are dominated by the same units. For $r_N>0$, let $\widecheck{\be}_i=\be_i/r_N$, so that $\sum_i\widecheck{\be}_i=\bzero$ by \eqref{eq:orth} and $\sum_i\|\widecheck{\be}_i\|^2=1$.

\begin{definition}
    \label{def:config}
    Let $q\ge1$, and for $N\ge2$ let $\bv_{N,1},\dots,\bv_{N,N}$ be vectors in $\R^q$ such that $\sup_N\sum_i\|\bv_{N,i}\|^2$ is finite. Let $\Psi=(\bpsi_k)_{k\in\mathcal{K}}$ be a finite or countably infinite family of nonzero vectors in $\R^q$ such that $\mathcal{K}(\zeta)=\{k:\|\bpsi_k\|>\zeta\}$ is finite for every $\zeta>0$, let $\mSig_E$ be a positive semidefinite $q\times q$ matrix, and let $\mathcal{I}_N(\zeta)=\{i:\|\bv_{N,i}\|>\zeta\}$. The array has \emph{limit configuration} $(\Psi,\mSig_E)$ if (i) and (ii) hold. (i) For every $\zeta>0$ different from all $\|\bpsi_k\|$ and all large $N$, $\mathcal{I}_N(\zeta)$ has the same cardinality $n(\zeta)$ as $\mathcal{K}(\zeta)$, and there are enumerations of the two sets, $i_1(N),\dots,i_{n(\zeta)}(N)$ and $k_1(N),\dots,k_{n(\zeta)}(N)$, with $\max_{j\le n(\zeta)}\|\bv_{N,i_j(N)}-\bpsi_{k_j(N)}\|\to0$. (ii) The limit $\lim_{\zeta\downarrow0}\limsup_{N\to\infty}\|\sum_{i\notin\mathcal{I}_N(\zeta)}\bv_{N,i}\bv_{N,i}^\top-\mSig_E\|=0$, the limit in $\zeta$ being taken along values different from all $\|\bpsi_k\|$.
\end{definition}

\begin{condition}[Residual configuration]
    \label{cond:E}
    For every $N$, $r_N>0$, and the array $(\widecheck{\be}_i)_{i\le N}$ has a limit configuration $(\Psi,\mSig_E)$ with $q=2$, where $\bpsi_k=(\psi_{k,1},\psi_{k,0})^\top$.
\end{condition}

The points $\bpsi_k$ are the limits of the normalized residual pairs of the units that keep a share of $r_N^2$ that does not vanish, and $\mSig_E$ is the limit of the second-moment matrix of the remaining units. If $\Psi$ is empty, $\max_i\|\widecheck{\be}_i\|\to0$, which is the negligibility condition of \citet{Noether1949} for the residual pairs. Let $\sigma_{E,1}^2$ and $\sigma_{E,0}^2$ be the diagonal entries of $\mSig_E$ corresponding to $e_i(1)$ and $e_i(0)$, let $\bw_\pi=(1/\pi,1/(1-\pi))^\top$, let $(A_k)_{k\in\mathcal{K}}$ be independent Bernoulli$(\pi)$ variables, and let $\bg_E\sim\mathcal{N}(\bzero,\pi(1-\pi)\mSig_E)$ be independent of $(A_k)$. Define
\begin{equation}
    \label{eq:UD}
    \begin{gathered}
        U_E=\sum_k(A_k-\pi)\bw_\pi^\top\bpsi_k+\bw_\pi^\top\bg_E,\\
        V_E=\sum_k\left\{\frac{A_k\psi_{k,1}^2}{\pi^2}+\frac{(1-A_k)\psi_{k,0}^2}{(1-\pi)^2}\right\}+\frac{\sigma_{E,1}^2}{\pi}+\frac{\sigma_{E,0}^2}{1-\pi}.
    \end{gathered}
\end{equation}

\begin{theorem}
    \label{thm:studentized}
    Suppose Assumption \ref{ass:design} and Conditions \ref{cond:lev} and \ref{cond:E} hold. Then,
    \begin{equation*}
        \left(\frac{N}{r_N}(\that-\tau_N),\ \frac{N^2}{r_N^2}\Vh\right)\wto(U_E,V_E).
    \end{equation*}
    Also, $\pr(V_E>0)=1$ if and only if $\mSig_E\ne\bzero$, or $\Psi$ is infinite, or some $\bpsi_k$ has both coordinates nonzero. In that case, $T_N\wto U_E/V_E^{1/2}$.
\end{theorem}

Section S4.1 proves Theorem \ref{thm:studentized}. If $\Psi$ is empty, $U_E$ is normal with mean zero and variance $\pi(1-\pi)\bw_\pi^\top\mSig_E\bw_\pi$, and $V_E$ is the constant $\sigma_{E,1}^2/\pi+\sigma_{E,0}^2/(1-\pi)$. The constant exceeds the variance of $U_E$ by $(1,-1)\mSig_E(1,-1)^\top\ge0$. Thus, $T_N$ has a normal limit with variance at most one, and the Wald interval is asymptotically conservative, as in \citet{Lin2013} and \citet{LiDing2017}. If $\Psi$ is not empty and $\pr(V_E>0)=1$, the limit of $T_N$ is nonnormal in general, and it depends on $\Psi$ as a family of points in $\R^2$, not only on its two coordinate projections. A unit whose treatment residual is large but whose control residual is not enters the numerator through $c_i$ in both arms, because $\tau_N$ contains both of its potential outcomes, but it enters the denominator only if it is treated. If $e_i(1)=e_i(0)$ for all $i$ and $\pi=1/2$, then $\bpsi_k=\alpha_k(1,1)^\top/2^{1/2}$ for real $\alpha_k$, and
\begin{equation}
    \label{eq:symmetric}
    \frac{U_E}{V_E^{1/2}}=\frac{\sum_k\alpha_k(2A_k-1)+\alpha_Eg}{\left(\sum_k\alpha_k^2+\alpha_E^2\right)^{1/2}},
\end{equation}
where $2A_k-1$ are independent Rademacher signs, $g\sim\mathcal{N}(0,1)$ is independent of them, and $\alpha_E=\{(1,1)\mSig_E(1,1)^\top/2\}^{1/2}$. The right-hand side of \eqref{eq:symmetric} is a sum of Rademacher signs with fixed weights plus an independent normal variable, and its variance is one. It is not normal if $\Psi$ is not empty, because its fourth cumulant is $-2\sum_k\alpha_k^4/(\sum_k\alpha_k^2+\alpha_E^2)^2<0$. Without equal residuals, the limit can be normal when $\Psi$ is not empty. If every $\bpsi_k$ is proportional to $(\pi,-(1-\pi))^\top$ and $\bw_\pi^\top\mSig_E\bw_\pi>0$, then the points cancel from $U_E$ and contribute only a constant to $V_E$.

\subsection{Heavy-tailed outcome errors}
\label{sec:sampling}

Configurations with infinitely many points arise when outcome errors are drawn from a heavy-tailed law. Let $\mathbb{S}^1$ be the unit circle in $\R^2$.

\begin{definition}
    \label{def:rv}
    Let $\gamma>0$. A random vector $\bvarepsilon\in\R^2$ is \emph{regularly varying with index $\gamma$ and spectral measure $\Lambda$}, a Borel probability measure on $\mathbb{S}^1$, if $\pr(\|\bvarepsilon\|>x)$ is regularly varying at infinity with index $-\gamma$ and $\pr(\bvarepsilon/\|\bvarepsilon\|\in\cdot|\|\bvarepsilon\|>x)$ converges weakly to $\Lambda$ as $x\to\infty$.
\end{definition}

Let $\mathsf{T}_1<\mathsf{T}_2<\cdots$ be the arrival times of a unit-rate Poisson process, let $\bso_1,\bso_2,\dots$ be independent with law $\Lambda$ and independent of $(\mathsf{T}_k)$, and define the random family $\Psi_\gamma=(\bpsi_k)_{k\ge1}$ by
\begin{equation}
    \label{eq:poissonpoints}
    \bpsi_k=\frac{\mathsf{T}_k^{-1/\gamma}\bso_k}{\left(\sum_{j\ge1}\mathsf{T}_j^{-2/\gamma}\right)^{1/2}} .
\end{equation}
Let $U_E$ and $V_E$ be defined by \eqref{eq:UD} with $\Psi=\Psi_\gamma$, $\mSig_E=\bzero$, and $(A_k)$ independent of $\Psi_\gamma$.

\begin{proposition}
    \label{prop:sampling}
    Let $(\bX_i,\bvarepsilon_i)$ for $i\ge1$ be independent and identically distributed, where $\bvarepsilon_i=(\varepsilon_i(1),\varepsilon_i(0))^\top$, and let $\Pi_N$ consist of the first $N$ units, with $Y_i(z)=\mu_z(\bX_i)+\varepsilon_i(z)$ for measurable functions $\mu_z$ such that $\E\mu_z(\bX_1)^2<\infty$, where $\mu_z$ is a constant if $d=0$. Suppose that $\bvarepsilon_1$ is independent of $\bX_1$, has mean zero, and is regularly varying with index $\gamma\in(1,2)$ and spectral measure $\Lambda$. If $d\ge1$, suppose also that $\E\|\bX_1\|^2<\infty$ and that the covariance matrix of $\bX_1$ is nonsingular. Let the assignment satisfy Assumption \ref{ass:design} and be independent of the population. On the event that $\mG$ is singular, let $h_i=0$ and $e_i(z)=Y_i(z)-\overline{Y}(z)$, and let ratios with denominator $r_N=0$ be zero.
    \begin{enumerate}
        \item[(a)] The share $\max_i\|\be_i\|^2/r_N^2$ converges in distribution to the law of $\mathsf{T}_1^{-2/\gamma}/\sum_j\mathsf{T}_j^{-2/\gamma}$, which is concentrated on $(0,1)$. If $d\ge1$, also $\pr(\mG\text{ singular})\to0$ and $\max_ih_i\to0$ in probability.
        \item[(b)] $\pr(V_E>0)=1$, and for every continuous function $\mathsf{f}$ on $[-\infty,\infty]$,
        \begin{equation}
            \label{eq:sampling}
            \E\{\mathsf{f}(T_N)|\Pi_N\}\wto\E\{\mathsf{f}(U_E/V_E^{1/2})|\Psi_\gamma\}.
        \end{equation}
        \item[(c)] The random variable $\E\{\mathsf{f}(U_E/V_E^{1/2})|\Psi_\gamma\}$ is almost surely constant for every continuous function $\mathsf{f}$ on $[-\infty,\infty]$ if and only if $\Lambda$ is concentrated on the line $\{\bm{x}\in\R^2:\bw_\pi^\top\bm{x}=0\}$.
    \end{enumerate}
\end{proposition}

Section S4.2 proves Proposition \ref{prop:sampling}, using the limits of the ranked norms and directions in \citet[Lemmas 1 and 3]{LePage1981}. In \eqref{eq:sampling}, $\E(\cdot|\Pi_N)$ is the randomization expectation with the population held fixed, the population being drawn once. By part (a), the largest residual pair keeps a share of $r_N^2$ that does not vanish, so the negligibility condition fails. By part (b), the randomization law of $T_N$ for the population at hand converges, in the sense of \eqref{eq:sampling}, to the law of Theorem \ref{thm:studentized} with $\mSig_E=\bzero$, whose points are the normalized points \eqref{eq:poissonpoints} of a Poisson process. This limit is a random law. By part (c), it varies with the population unless $\Lambda$ is concentrated on the line $\bw_\pi^\top\bm{x}=0$. In that case, every point lies on this line, the points cancel from $U_E$, and the limit is the point mass at zero. The limits of the bootstrap of the mean with infinite variance \citep{Athreya1987} and of permuted sums of heavy-tailed samples \citep{AueBerkesHorvath2008, BerkesHorvathSchauer2010} also depend on the sample. The functions $\mu_z$ allow conditional means that are nonlinear in $\bX_i$, and their projection residuals enter $r_N^2$ only at the order $N$, which is smaller than the order of the heavy-tailed errors.

\subsection{Pairing of the potential outcomes}
\label{sec:pairing}

Each unit reveals one potential outcome, so the pairing of $Y_i(1)$ with $Y_i(0)$ is unobserved, and the variance of the unit-level effects is thus unidentified \citep{AronowGreenLee2014}. The next proposition shows that this pairing can change the limit of $T_N$.

\begin{proposition}
    \label{prop:pairing}
    Let $d=0$, suppose Assumption \ref{ass:design} holds with $N_1=N/2$ along even $N$, and let $x_N>0$ satisfy $x_N/N^{1/2}\to\infty$, and let $w_i=(-1)^i$. In population I, let $Y_1(1)=Y_1(0)=x_N$, $Y_2(1)=Y_2(0)=0$ and $Y_i(1)=Y_i(0)=w_i$ for $i\ge3$. In population II, let $Y_1(1)=x_N$, $Y_1(0)=0$, $Y_2(1)=0$, $Y_2(0)=x_N$ and $Y_i(1)=Y_i(0)=w_i$ for $i\ge3$. Then, both populations have $\tau_N=0$ and the same multisets $\{Y_i(1):i\le N\}$ and $\{Y_i(0):i\le N\}$, but, along even $N$,
    \begin{equation*}
        \begin{gathered}
            \pr_{\mathrm{I}}(|T_N|\le t)\to0,
            \qquad
            \pr_{\mathrm{II}}(|T_N|\le t)\to\Phi(t)/2
            \qquad
            (0<t<1),\\
            \pr_{\mathrm{I}}(|T_N|\le t)\to1,
            \qquad
            \pr_{\mathrm{II}}(|T_N|\le t)\to\{1+\Phi(t)\}/2
            \qquad
            (t>1),
        \end{gathered}
    \end{equation*}
    where $\pr_{\mathrm{I}}$ and $\pr_{\mathrm{II}}$ are the randomization distributions for the two populations.
\end{proposition}

Section S4.4 proves Proposition \ref{prop:pairing}. In the regular case, the pairing affects only $(1,-1)\mSig_E(1,-1)^\top$, the limit of $\sum_i\{e_i(1)-e_i(0)\}^2/r_N^2$. The limit of the normalized variance estimator exceeds the variance of the limit of the normalized estimator by this amount, so the Wald interval is conservative. With a residual configuration, it can instead change the shape of the limit law. The two populations of Proposition \ref{prop:pairing} have different limiting coverage at every critical value $t\ne1$, and the ordering of the two coverages reverses at $t=1$. Population I has the single point $(1,1)^\top/2^{1/2}$ and $\mSig_E=\bzero$. Population II has the points $(1,0)^\top/2^{1/2}$ and $(0,1)^\top/2^{1/2}$ and $\mSig_E=\bzero$, so $\pr(V_E=0)=1/4$ in Theorem \ref{thm:studentized}.

\section{Coverage under balanced allocation}
\label{sec:floor}

The limits of Theorem \ref{thm:studentized} depend on how $e_i(1)$ is paired with $e_i(0)$ at the dominant units, which cannot be identified from the observed data. Under balanced allocation, that is, under Assumption \ref{ass:design} with $\pi=1/2$, the coverage of the Wald interval nevertheless admits explicit lower bounds. Let $\Phibar=1-\Phi$.

\begin{proposition}
    \label{prop:balancedfloor}
    Suppose Assumption \ref{ass:design} holds with $\pi=1/2$, and Condition \ref{cond:lev} holds. Then, for every $t>0$,
    \begin{equation}
        \label{eq:balancedfloor}
        \liminf_{N\to\infty}\pr(|T_N|\le t)\ge\frac{1-2\cstar\Phibar(t)}{2},
        \qquad
        \cstar=\frac{1}{4\Phibar(2^{1/2})}\approx3.1787.
    \end{equation}
    In particular, the Wald interval of nominal level $1-\varsigma$ has asymptotic coverage at least $(1-\cstar\varsigma)/2$, which is approximately $0.421$ at $\varsigma=0.05$.
\end{proposition}

Section S5.2 proves Proposition \ref{prop:balancedfloor}. The bound requires no condition on the outcomes. When $N_1=N/2$, the assignments $(Z_1,\dots,Z_N)$ and $(1-Z_1,\dots,1-Z_N)$ are equally likely. Thus, when $\pi_N\to1/2$, changing the second assignment at a vanishing fraction of the units gives an assignment with the law of the first. At the two assignments, the numerators of $T_N$ have the same absolute value to first order, and the sums of the squared population residuals over the observed arms add up to $r_N^2$ to first order. For one of the two assignments, the variance estimate is therefore at least about $2r_N^2/N^2$, and the interval covers $\tau_N$ unless a weighted sum of the signs $2Z_i-1$, with squared weights summing to at most one, is large. Up to a term that vanishes in probability, this sum is a sum of independent Rademacher signs with the same weights, and $\cstar$ is the sharp constant of the Gaussian comparison for such sums of \citet{BentkusDzindzalieta2015}.

A higher floor holds if the two potential residuals of the units are matched in aggregate. Let $m_i=e_i(1)+e_i(0)$.

\begin{condition}[Aggregate matching]
    \label{cond:F}
    There is $\eta>0$ such that
    \begin{equation*}
        \liminf_N\frac{\sum_i\min\{e_i(1)^2,e_i(0)^2\}}{\sum_im_i^2}\ge\eta^2,
    \end{equation*}
    where the ratio is $+\infty$ if $\sum_im_i^2=0$.
\end{condition}

If $e_i(1)=e_i(0)$ for all $i$, Condition \ref{cond:F} holds with $\eta=1/2$. If $e_i(1)=\rho_ie_i(0)$ with $\rho_i\in[\underline{\rho},1/\underline{\rho}]$ for some $\underline{\rho}\in(0,1]$, it holds with $\eta=\underline{\rho}/(1+\underline{\rho})$ (Section S5.3). Geometrically, plot each unit as the point $(e_i(1),e_i(0))$. Under balanced allocation, each unit enters the numerator of $T_N$ through $m_i$, which is proportional to the projection of its point onto the diagonal, with a sign given by its arm. It enters the variance estimator through the coordinate of its arm, whose square is at least $\min\{e_i(1)^2,e_i(0)^2\}$, the squared distance from the point to the nearer coordinate axis. Condition \ref{cond:F} prevents the residual mass from concentrating, in aggregate, near the axes relative to its projection onto the diagonal. Mass near the line $e_i(1)=-e_i(0)$ enters the variance estimator but contributes little to the numerator, and it increases the ratio in Condition \ref{cond:F}.

\begin{theorem}
    \label{thm:floor}
    Suppose Assumption \ref{ass:design} holds with $\pi=1/2$, and Conditions \ref{cond:lev} and \ref{cond:F} hold. Then, for every $t>0$,
    \begin{equation}
        \label{eq:floor}
        \limsup_{N\to\infty}\pr(|T_N|>t)\le2\cstar\Phibar(2\eta t).
    \end{equation}
    In particular, if $e_i(1)=e_i(0)$ for all $i$, the Wald interval of nominal level $1-\varsigma$ has asymptotic coverage at least $1-\cstar\varsigma$.
\end{theorem}

Section S5.3 proves Theorem \ref{thm:floor}. Theorem \ref{thm:floor} requires no configuration condition on the residuals. For $\eta=1/2$, the floor is approximately $0.841$ at $\varsigma=0.05$ and $0.968$ at $\varsigma=0.01$. For $\eta=1/3$, which corresponds to pointwise matching with $\underline{\rho}=1/2$, it is approximately $0.392$ at $\varsigma=0.05$. At $\varsigma=0.05$, the floor of Theorem \ref{thm:floor} exceeds that of Proposition \ref{prop:balancedfloor} only if $\eta>0.340$. The proof replaces complete randomization by independent signs through a coupling of the type of \citet{Hajek1960}. The magnitudes $|m_i|$ are then fixed, as the absolute values of the observations are in the analysis of Student's $t$ statistic under symmetry \citep{Efron1969}, and the bound has the form of the tail bound of \citet{Pinelis1994} for that statistic under symmetry, with the constant of \citet{BentkusDzindzalieta2015}.

For $K\ge1$ and $0\le j\le K$, let $\theta_K(j)=\{j/\pi-(K-j)/(1-\pi)\}/\{j/\pi^2+(K-j)/(1-\pi)^2\}^{1/2}$, where $\pi$ is the limit in Assumption \ref{ass:design}. The next proposition shows that both floors are close to the best possible and that neither holds at the allocation $\pi=0.4$ and the nominal level $0.95$.

\begin{proposition}
    \label{prop:counter}
    Let $d=0$, suppose Assumption \ref{ass:design} holds, let $K\ge1$ be fixed, let $x_N>0$ satisfy $x_N/N^{1/2}\to\infty$, and let $w_{N,i}$ satisfy $\sup_NN^{-1}\sum_{i\le N}w_{N,i}^2<\infty$.
    \begin{enumerate}
        \item[(a)] If $Y_i(1)=Y_i(0)=x_N$ for $i\le K$ and $Y_i(1)=Y_i(0)=w_{N,i}$ for $i>K$, then for every $t>0$ outside the finite set $\{|\theta_K(j)|:0\le j\le K\}$,
        \begin{equation}
            \label{eq:counterA}
            \pr(|T_N|\le t)\to\sum_{j=0}^K\binom Kj\pi^j(1-\pi)^{K-j}\1\{|\theta_K(j)|\le t\}.
        \end{equation}
        For $K=4$ and $t=\Phi^{-1}(0.975)$, $\pr(|T_N|\le t)\to1-\pi^4-(1-\pi)^4$ for every $\pi\in(0,1)$. The limit is $7/8$ at $\pi=1/2$ and approximately $0.589$ at $\pi=0.2$.
        \item[(b)] Let $Y_i(1)=Y_i(0)=w_{N,i}$ for $i\ge2$. If $Y_1(1)=x_N$ and $Y_1(0)=0$, then $\pr(|T_N|\le t)\to\pi$ for every $t>1-\pi$ and $\pr(|T_N|\le t)\to0$ for every $t\in(0,1-\pi)$. If $Y_1(1)=0$ and $Y_1(0)=x_N$, the same holds with $\pi$ and $1-\pi$ interchanged.
        \item[(c)] Let $\pi=0.4$, $Y_1(1)=Y_2(1)=5x_N$, $Y_1(0)=-2x_N$, $Y_2(0)=-x_N$, and $Y_i(1)=Y_i(0)=w_{N,i}$ for $i\ge3$. Then, Condition \ref{cond:F} holds with $\eta=5^{-1/2}$, and for $t=\Phi^{-1}(0.975)$, $\pr(|T_N|\le t)\to1-(1-\pi)^2=0.64$.
    \end{enumerate}
\end{proposition}

Section S5.4 proves Proposition \ref{prop:counter}. In part (a), Condition \ref{cond:F} holds with $\eta=1/2$, and $|\theta_K(0)|=|\theta_K(K)|=K^{1/2}$, whereas $|\theta_K(j)|<(K-1)^{1/2}$ for $0<j<K$. For $K=4$ and $t=\Phi^{-1}(0.975)$, the assignments that place all four dominant units in one arm are therefore not covered, and all other assignments are. At $\pi=1/2$, the limiting coverage is $7/8$, so no floor valid for every population with $\eta=1/2$ exceeds $7/8$ at $\varsigma=0.05$, and the floor $0.841$ of Theorem \ref{thm:floor} is close to sharp. The probability $\pi^4+(1-\pi)^4$ of these assignments tends to one as $\pi\to0$, and thus aggregate matching without balanced allocation gives no floor that holds for every $\pi$. With $K=2$ and $\pi=1/2$, $|\theta_2(j)|\in\{0,2^{1/2}\}$, so the limiting coverage is $1/2$ for $0<t<2^{1/2}$, whereas the floor $1-2\cstar\Phibar(t)$ of Theorem \ref{thm:floor} with $\eta=1/2$ tends to $1-2\cstar\Phibar(2^{1/2})=1/2$ as $t\uparrow2^{1/2}$. The constant $2\cstar$ in \eqref{eq:floor} therefore cannot be lowered.

In part (b), the single dominant unit lies on a coordinate axis, so the ratio in Condition \ref{cond:F} tends to zero. If that unit is treated, $T_N\to1-\pi$ in probability and the interval covers. If it is assigned to control, its observed outcome $Y_1(0)=0$ is of the size of the outcomes of the remaining units, so $N^2\Vh=O(N)$, whereas $N(\that-\tau_N)=-x_N+O_p(N^{1/2})$, and $|T_N|\to\infty$. At $\pi=1/2$, no floor valid for every population under balanced allocation exceeds $1/2$, and the universal floor $0.421$ of Proposition \ref{prop:balancedfloor} is within a factor of $1.19$ of this bound. For general $\pi$, the two populations of part (b) have limiting coverage $\pi$ and $1-\pi$, so without balanced allocation the coverage can be as low as $\min(\pi,1-\pi)$. In part (c), only the assignments that place both dominant units in control are not covered. The floor of Theorem \ref{thm:floor} at $\eta=5^{-1/2}$ is approximately $0.747$, above the limiting coverage $0.64$, and for the first population of part (b) at $\pi=0.4$ the limiting coverage $0.4$ is below the floor $0.421$ of Proposition \ref{prop:balancedfloor}. Hence neither floor holds at the allocation $\pi=0.4$ and the nominal level $0.95$.

\section{Numerical studies}
\label{sec:sim}

\subsection{Finite populations and the limit laws}
\label{sec:simlimits}

This study checks how closely the limits of Theorem \ref{thm:studentized} and Propositions \ref{prop:pairing} and \ref{prop:counter} describe finite populations of moderate size. The populations of Propositions \ref{prop:pairing} and \ref{prop:counter} have $d=0$, $N=20{,}000$ and $x_N=N$, with $w_{N,i}=(-1)^i$, $K=4$ in Proposition \ref{prop:counter}(a), and $\pi\in\{0.5,0.2\}$ in Proposition \ref{prop:counter}. For Theorem \ref{thm:studentized}, the population has $N=20{,}000$, $\pi=0.3$, one standard normal covariate, and potential outcomes $Y_i(1)=2X_i+\varepsilon_i(1)$ and $Y_i(0)=X_i+\varepsilon_i(0)$. The error pairs are bivariate normal with unit variances and correlation $0.5$, except for three units, whose error pairs are $N^{1/2}(2,0.5)$, $N^{1/2}(-1.5,1)$ and $N^{1/2}(0.5,-2)$. Its limit law is evaluated at its own three largest normalized residual pairs and the second-moment matrix of the remaining pairs, from $10^6$ draws. For each population, the study draws $100{,}000$ assignments and records the coverage of $\that\pm t\Vh^{1/2}$ for $t\in\{0.25,1.96\}$. Table \ref{tab:limits} reports these coverages next to their limits.

\begin{table}[ht]
    \centering
    \caption{Coverage of $\that\pm t\Vh^{1/2}$ over $100{,}000$ assignments in populations of size $N=20{,}000$, and its limit. The columns give the allocation $\pi$, the largest share $\max_i\|\be_i\|^2/r_N^2$ (share), and, for $t=0.25$ and $t=1.96$, the coverage over assignments (cov.) and the limit of the cited result (limit). The largest Monte Carlo standard error is 0.0016.}
    \label{tab:limits}
    \begin{tabular}{lcccccc}
        \toprule
        {} & {} & {} & \multicolumn{2}{c}{$t=0.25$} & \multicolumn{2}{c}{$t=1.96$}\\
        \cmidrule(lr){4-5}\cmidrule(lr){6-7}
        population & $\pi$ & share & cov. & limit & cov. & limit\\
        \midrule
        Proposition \ref{prop:pairing}, population I & 0.5 & 1.000 & 0.000 & 0.000 & 1.000 & 1.000\\
        Proposition \ref{prop:pairing}, population II & 0.5 & 0.500 & 0.301 & 0.299 & 0.987 & 0.988\\
        Proposition \ref{prop:counter}(a), $K=4$ & 0.5 & 0.250 & 0.372 & 0.375 & 0.875 & 0.875\\
        Proposition \ref{prop:counter}(a), $K=4$ & 0.2 & 0.250 & 0.407 & 0.410 & 0.588 & 0.589\\
        Proposition \ref{prop:counter}(b) & 0.5 & 1.000 & 0.000 & 0.000 & 0.502 & 0.500\\
        Proposition \ref{prop:counter}(b) & 0.2 & 1.000 & 0.000 & 0.000 & 0.198 & 0.200\\
        Theorem \ref{thm:studentized}, $d=1$ & 0.3 & 0.309 & 0.203 & 0.204 & 0.998 & 0.998\\
        \bottomrule
    \end{tabular}
\end{table}

The results match Theorem \ref{thm:studentized} and Propositions \ref{prop:pairing} and \ref{prop:counter}. Every coverage agrees with its limit within two Monte Carlo standard errors, including the limits zero and one. The two populations of Proposition \ref{prop:pairing} have the same marginal distributions of potential outcomes, yet their coverages at both critical values differ as the limits predict. In the population of Proposition \ref{prop:counter}(b), the coverage of the nominal $95\%$ interval is close to the allocation $\pi$.

\subsection{Coverage with heavy-tailed errors}
\label{sec:simcoverage}

This study examines the coverage of the nominal $95\%$ Wald interval when the errors are heavy-tailed. Each population has $N=4{,}000$ units with three independent standard normal covariates, $Y_i(1)=\bX_i^\top(2,2,2)^\top+\varepsilon_i(1)$ and $Y_i(0)=\bX_i^\top(1,1,1)^\top+\varepsilon_i(0)$. A symmetric Pareto variable with index $\gamma$ is generated as $\mathsf{P}\mathsf{s}$, where $\pr(\mathsf{P}>x)=x^{-\gamma}$ for $x\ge1$ and $\mathsf{s}$ is an independent random sign. The errors are symmetric Pareto with index $\gamma\in\{1.2,1.5,1.8\}$ and of three kinds: (i) additive errors that have $\varepsilon_i(1)=\varepsilon_i(0)$; (ii) matched errors that have $\varepsilon_i(1)=\rho_i\varepsilon_i(0)$, with $\log\rho_i$ uniform on $[\log(1/2),\log2]$; and (iii) independent errors that have $\varepsilon_i(1)$ and $\varepsilon_i(0)$ independent. The allocation is $\pi\in\{0.5,0.2\}$. Each setting has $200$ populations, and each population has $1{,}000$ assignments, analyzed with the interacted estimator on the three covariates and with the difference in means. For each population, the study computes $\max_i\|\be_i\|^2/r_N^2$, $\eta_N^2=\sum_i\min\{e_i(1)^2,e_i(0)^2\}/\sum_im_i^2$, which is the ratio of Condition \ref{cond:F} at $N$, the floor $1-2\cstar\Phibar\{2\eta_N\Phi^{-1}(0.975)\}$ of Theorem \ref{thm:floor}, and the coverage. Table \ref{tab:coverage} reports summaries of these quantities over the populations of each setting.

The results match Propositions \ref{prop:balancedfloor} and \ref{prop:counter} and Theorem \ref{thm:floor}. Under balanced allocation, no population has coverage more than two Monte Carlo standard errors below the larger of its floor of Theorem \ref{thm:floor} and the floor $0.421$ of Proposition \ref{prop:balancedfloor}. The largest residual pair carries a large share of $r_N^2$, and this share decreases as $\gamma$ grows. With additive and matched errors, the mean coverage is at or above $0.95$, and the lowest coverage stays close to it under balanced allocation. With independent errors, the median of $\eta_N$ is small and the floor of Theorem \ref{thm:floor} at this median is not positive. Specifically, independent heavy-tailed errors place the dominant residual pairs near the coordinate axes, so $\eta_N$ falls below $0.257$, where the floor of Theorem \ref{thm:floor} is negative and only the floor $0.421$ of Proposition \ref{prop:balancedfloor} is informative. The mean coverage is below $0.95$ for $\gamma=1.2$, and for $\gamma=1.5$ at $\pi=0.2$. The lowest coverage under balanced allocation is then close to $1/2$, above the floor of Proposition \ref{prop:balancedfloor}, and the populations with the lowest coverage have one residual pair that carries almost all of $r_N^2$, as in Proposition \ref{prop:counter}(b). At $\pi=0.2$, no floor applies, and with independent errors and $\gamma\le1.5$ the lowest coverage falls close to $\pi$. The interacted estimator and the difference in means have nearly the same mean coverage in every setting, because the covariates contribute to $r_N^2$ only at the order $N$, which is smaller than the order of the heavy-tailed errors (Section \ref{sec:sampling}).

\begin{table}[t]
    \centering
    \caption{Coverage of nominal $95\%$ Wald intervals with symmetric Pareto errors of index $\gamma$, $N=4{,}000$, and $200$ populations of $1{,}000$ assignments each, for the interacted estimator on three normal covariates ($\that$) and the difference in means ($d=0$). The columns give, for the interacted estimator, the median of $\eta_N$ over populations, the floor of Theorem \ref{thm:floor} at this median and the median of $\max_i\|\be_i\|^2/r_N^2$ (share), and for each estimator the mean and the minimum of the coverage and the number of populations whose coverage is more than two Monte Carlo standard errors below the larger of their own floor of Theorem \ref{thm:floor} and the floor $0.421$ of Proposition \ref{prop:balancedfloor} (below). The Monte Carlo standard error at coverage $x$ is $\{x(1-x)/1000\}^{1/2}$, evaluated at that larger floor. No floor applies for $\pi=0.2$ (--). The largest leverage in any population is 0.0097.}
    \label{tab:coverage}
    \setlength{\tabcolsep}{3.5pt}
    \begin{tabular}{lllccccccccc}
        \toprule
        {} & {} & {} & {} & {} & {} & \multicolumn{3}{c}{$\that$} & \multicolumn{3}{c}{$d=0$}\\
        \cmidrule(lr){7-9}\cmidrule(lr){10-12}
        $\gamma$ & errors & $\pi$ & $\eta_N$ & floor & share & mean & min & below & mean & min & below\\
        \midrule
        1.2 & additive & 0.5 & 0.500 & 0.841 & 0.516 & 0.973 & 0.933 & 0 & 0.972 & 0.932 & 0\\
        1.2 & additive & 0.2 & 0.500 & -- & 0.466 & 0.972 & 0.907 & -- & 0.972 & 0.914 & --\\
        1.2 & matched & 0.5 & 0.414 & 0.668 & 0.508 & 0.967 & 0.934 & 0 & 0.967 & 0.933 & 0\\
        1.2 & matched & 0.2 & 0.415 & -- & 0.515 & 0.960 & 0.685 & -- & 0.960 & 0.690 & --\\
        1.2 & independent & 0.5 & 0.040 & $<0$ & 0.511 & 0.886 & 0.483 & 0 & 0.886 & 0.483 & 0\\
        1.2 & independent & 0.2 & 0.045 & -- & 0.490 & 0.879 & 0.180 & -- & 0.880 & 0.180 & --\\
        1.5 & additive & 0.5 & 0.500 & 0.841 & 0.378 & 0.964 & 0.938 & 0 & 0.963 & 0.941 & 0\\
        1.5 & additive & 0.2 & 0.500 & -- & 0.356 & 0.963 & 0.915 & -- & 0.962 & 0.929 & --\\
        1.5 & matched & 0.5 & 0.418 & 0.678 & 0.350 & 0.964 & 0.937 & 0 & 0.964 & 0.938 & 0\\
        1.5 & matched & 0.2 & 0.413 & -- & 0.329 & 0.956 & 0.839 & -- & 0.957 & 0.845 & --\\
        1.5 & independent & 0.5 & 0.126 & $<0$ & 0.342 & 0.954 & 0.484 & 0 & 0.955 & 0.486 & 0\\
        1.5 & independent & 0.2 & 0.134 & -- & 0.288 & 0.933 & 0.202 & -- & 0.934 & 0.202 & --\\
        1.8 & additive & 0.5 & 0.500 & 0.841 & 0.198 & 0.955 & 0.934 & 0 & 0.957 & 0.932 & 0\\
        1.8 & additive & 0.2 & 0.500 & -- & 0.209 & 0.956 & 0.922 & -- & 0.956 & 0.930 & --\\
        1.8 & matched & 0.5 & 0.412 & 0.663 & 0.199 & 0.958 & 0.940 & 0 & 0.959 & 0.938 & 0\\
        1.8 & matched & 0.2 & 0.413 & -- & 0.251 & 0.955 & 0.867 & -- & 0.956 & 0.895 & --\\
        1.8 & independent & 0.5 & 0.235 & $<0$ & 0.187 & 0.970 & 0.484 & 0 & 0.969 & 0.491 & 0\\
        1.8 & independent & 0.2 & 0.232 & -- & 0.214 & 0.974 & 0.531 & -- & 0.969 & 0.691 & --\\
        \bottomrule
    \end{tabular}
\end{table}

\section{Data application}
\label{sec:data}

This section applies the results to data from three randomized experiments whose outcomes, medical spending, emergency-department charges, and customer spending, are heavy-tailed. The coverage of an interval cannot be evaluated from the experiment itself, because only one assignment is observed and $\tau_N$ is unknown. We therefore impute the potential outcomes of each experiment under a sharp model of the effect and draw new assignments by complete randomization. The coverage is then the proportion of assignments whose interval contains $\tau_N$.

The first population comes from the RAND Health Insurance Experiment \citep{AronDineEinavFinkelstein2013}. We use the person-year file distributed with \citet{AngristPischke2015} and keep the $5{,}477$ persons observed in both of their first two years in the experiment. The treated arm consists of the $1{,}843$ persons assigned to free care, and the control arm consists of those assigned to a plan with cost sharing. The outcome $Y_i$ is the person's total medical spending in the second year, and $X_i$ is the total in the first year. The second population comes from the Oregon Health Insurance Experiment \citep{FinkelsteinEtAl2012, TaubmanEtAl2014}. We use the public-use emergency-department file for the Portland area and keep the $19{,}619$ lottery participants who signed up alone and have both charge totals recorded. The treated arm consists of the $6{,}750$ persons selected in the lottery. The outcome $Y_i$ is the sum of the charges for the person's emergency-department visits from 10 March 2008 to 30 September 2009, and $X_i$ is the same sum from 1 January 2007 to 9 March 2008, before the first lottery notification. The third population comes from the e-mail experiment of \citet{Hillstrom2008}, in which $64{,}000$ customers of a retailer were assigned at random to a men's-merchandise e-mail, a women's-merchandise e-mail or no e-mail. We keep the $42{,}613$ customers assigned to the men's e-mail, which forms the treated arm, or to no e-mail. The outcome $Y_i$ is the customer's spending in the two weeks after the e-mail, and $X_i$ is the spending in the preceding year.

The recorded outcome of each unit is its potential outcome under the arm to which it was assigned in the experiment, and the other potential outcome is imputed under one of two sharp models. The constant effect sets $Y_i(1)-Y_i(0)=\tau_N$ for every unit, with $\tau_N$ equal to the difference between the average recorded outcomes of the treated and control arms. It gives $e_i(1)=e_i(0)$, so Condition \ref{cond:F} holds with $\eta=1/2$. The proportional effect sets $Y_i(1)=\rho\,Y_i(0)$ for every unit, with $\rho$ equal to the ratio of these averages. It gives $e_i(1)=\rho\,e_i(0)$, and since $\rho>1$ in the three experiments, Condition \ref{cond:F} holds with $\eta=1/(1+\rho)$. On the completed population, we draw $100{,}000$ assignments by complete randomization with $\pi=0.5$ and with $\pi=0.2$, and compute the nominal $95\%$ Wald interval for the difference in means ($d=0$) and for the interacted estimator. The interacted estimator uses the covariate $\bX_i=\log(1+X_i)$, so $d=1$. We do not use $X_i$ itself, because its leverage is about $0.19$ for one person in the RAND population, outside the scope of Condition \ref{cond:lev}. Because the assignments are redrawn, the covariate need not precede the original assignment, as in the RAND population. Table \ref{tab:data} reports the coverage, together with the maximum leverage of $\bX_i$, the largest share $\max_i\|\be_i\|^2/r_N^2$, the ratio $\eta_N$ of Condition \ref{cond:F} at $N$ defined in Section \ref{sec:simcoverage}, and the floor of Theorem \ref{thm:floor} at $\eta_N$.

In the RAND population, the largest residual pair carries a large share of $r_N^2$ for both estimators, whereas in the Oregon and retail populations, its share is small. The maximum leverage of $\bX_i$ is small in all three populations. Under balanced allocation, the coverage is above $0.95$ in every population, well above the floors of Proposition \ref{prop:balancedfloor} and Theorem \ref{thm:floor}. With $\pi=0.2$, the floors do not apply, and by Proposition \ref{prop:counter} they need not hold. The coverage in the RAND population then falls below $0.95$, more under the proportional effect than under the constant effect, whereas in the Oregon and retail populations, it stays close to $0.95$. The covariate changes the largest share only slightly, and the two estimators have nearly the same coverage in every setting. Thus, in these populations, the largest loss of coverage occurs under unbalanced allocation in the population in which a few units dominate the residual variation.

\begin{table}[t]
    \centering
    \caption{Coverage of nominal $95\%$ Wald intervals over $100{,}000$ assignments in the RAND ($N=5{,}477$), Oregon ($N=19{,}619$) and retail ($N=42{,}613$) populations under a constant or proportional (prop.) effect, for the difference in means ($d=0$) and for the interacted estimator on $\bX_i=\log(1+X_i)$ (log). The constant effect $\tau_N$ is $155.55$, $133.89$ and $0.77$, and the ratio $\rho$ of the proportional effect is $1.34$, $1.03$ and $2.18$. The number of treated units $N_1$ is $2{,}738$, $9{,}810$ and $21{,}306$ at $\pi=0.5$ and $1{,}095$, $3{,}924$ and $8{,}523$ at $\pi=0.2$. All values are given for the three populations in the order listed. The columns give the maximum leverage $\max_ih_i$ of $\bX_i$, the largest share $\max_i\|\be_i\|^2/r_N^2$ for each estimator, $\eta_N$, the floor of Theorem \ref{thm:floor} at $\eta_N$, which applies at $\pi=0.5$, and the coverage. The largest Monte Carlo standard error is 0.0008.}
    \label{tab:data}
    \setlength{\tabcolsep}{4pt}
    \begin{tabular}{llccccccccc}
        \toprule
        {} & {} & {} & \multicolumn{2}{c}{share} & {} & {} & \multicolumn{2}{c}{$\pi=0.5$} & \multicolumn{2}{c}{$\pi=0.2$}\\
        \cmidrule(lr){4-5}\cmidrule(lr){8-9}\cmidrule(lr){10-11}
        population & effect & $\max_ih_i$ & $d=0$ & log & $\eta_N$ & floor & $d=0$ & log & $d=0$ & log\\
        \midrule
        RAND & constant & 0.0014 & 0.329 & 0.338 & 0.500 & 0.841 & 0.954 & 0.954 & 0.937 & 0.937\\
        RAND & prop. & 0.0014 & 0.373 & 0.381 & 0.428 & 0.703 & 0.955 & 0.955 & 0.926 & 0.925\\
        Oregon & constant & 0.0003 & 0.016 & 0.016 & 0.500 & 0.841 & 0.951 & 0.951 & 0.947 & 0.947\\
        Oregon & prop. & 0.0003 & 0.016 & 0.016 & 0.492 & 0.829 & 0.951 & 0.951 & 0.947 & 0.947\\
        retail & constant & 0.0002 & 0.026 & 0.026 & 0.500 & 0.841 & 0.951 & 0.951 & 0.945 & 0.945\\
        retail & prop. & 0.0002 & 0.058 & 0.058 & 0.315 & 0.308 & 0.961 & 0.961 & 0.944 & 0.944\\
        \bottomrule
    \end{tabular}
\end{table}

\section{Discussion}
\label{sec:disc}

We studied the studentized difference in means and the interacted estimator, with covariates whose leverage vanishes, when a few units carry a share of the residual variation that does not vanish. The studentized statistic then has a limit that is nonnormal in general and depends on how the two potential residuals of the dominant units are paired. When the errors are heavy-tailed, this limit is random and, outside a degenerate case, varies with the population. Under balanced allocation, the Wald interval keeps explicit coverage floors, which are close to the best possible, and neither floor extends to every allocation.

Natural next steps are the best universal floor under balanced allocation, which lies between $0.421$ and the coverage $1/2$ of Proposition \ref{prop:counter}(b), a coverage floor for a fixed allocation $\pi\ne1/2$, and an interval whose coverage is guaranteed for every population, including populations whose covariates have leverage that does not vanish.

\section*{Data availability}

The RAND data are distributed with \citet{AngristPischke2015}, the Oregon data are the public-use files of the Oregon Health Insurance Experiment \citep{FinkelsteinEtAl2012}, and the retail data are those of \citet{Hillstrom2008}.

\section*{Supplementary Material}

Section S1 of the Supplementary Material states conventions and external results. Sections S2--S5 prove the results of Sections \ref{sec:setup}--\ref{sec:floor}. The Supplementary Material is provided as an ancillary file with this article.

\bibliographystyle{jasa3}
\bibliography{p2}

\end{document}